\documentclass[12pt]{article}
\baselineskip
18pt
\begin{document}
\baselineskip 24pt\baselineskip 24pt
\begin{center}

%\vspace{3in}

\Large
{\bf Generalized hyperbolic functions, circulant matrices
and functional equations}\\

\normalsize
\vspace{.5cm}
by\\
\bigskip
\bigskip
Martin E. Muldoon\footnote{Research supported by grants from the
Natural Sciences and Engineering Research Council, Canada }\\
Department of Mathematics \& Statistics\\
York University\\
Toronto ON M3J 1P3, Canada\\
\bigskip
email: {\tt muldoon@yorku.ca}\\
\bigskip
Revised March 4, 2005
\bigskip
\end{center}

%\newpage
\begin{abstract}
There is a contrast between the two sets of functional equations
$$\begin{array}{c}
f_0(x+y) = f_0(x)f_0(y) +  f_1(x)f_1(y), \\
f_1(x+y) = f_1(x)f_0(y) +  f_0(x)f_1(y) ,
\end{array}
$$ and
$$ \begin{array}{c}
f_0(x-y) = f_0(x)f_0(y) -  f_1(x)f_1(y),\\
f_1(x-y) = f_1(x)f_0(y) -  f_0(x)f_1(y)
\end{array} $$
satisfied by the even and odd components of a solution of $f(x+y) =
f(x) f(y)$. J. Schwaiger and, later, W. F\"org-Rob and J. Schwaiger
 considered the extension of these ideas to the case where $f$
is sum of $n$ components. Here we shorten and simplify the
statements and proofs of some of these
 results by a more systematic use of
matrix notation.
\end{abstract}
\bigskip
\bigskip



{\em Abbreviated title}: Generalized hyperbolic functions

\newpage
\section{Introduction}
\setcounter{equation}{0}
By a more systematic use of matrix notation, we shorten the
statements and proofs of some results of J. Schwaiger \cite{schw}
on  generalized
hyperbolic functions and their characterization by functional
equations.  We also discuss some stability results for these
equations; our results here are motivated by those of W.
F\"org--Rob and J. Schwaiger \cite{forg}; see also \cite{forg1}.

The results discussed here can be readily understood in the case
$n=2$. A function $f: C \rightarrow C$ can be written as a sum of
even and
odd components in a standard way:
$f(x) = f_0(x) + f_1(x), $
where
$f_0(x) = [f(x) + f(-x)]/2,\;\;f_1(x) =
[f(x) - f(-x)]/2.
$
In case $f$ is an {\em  exponential function}, i.e, if $f$ satisfies
\begin{equation} f(x+y) = f(x) f(y), \label{gef} \end{equation}
for all $x$ and $y$, then the  components satisfy
\begin{equation}
\begin{array}{c}
f_0(x+y) = f_0(x)f_0(y) +  f_1(x)f_1(y), \\
f_1(x+y) = f_1(x)f_0(y) +  f_0(x)f_1(y) ,
\end{array} \label{fplus}
\end{equation}
and
\begin{equation}
\begin{array}{c}
f_0(x-y) = f_0(x)f_0(y) -  f_1(x)f_1(y),\\
f_1(x-y) = f_1(x)f_0(y) -  f_0(x)f_1(y).
\end{array}
\label{fminus}
\end{equation}
In case $f$ is the usual exponential function, then $f_0$ and
$f_1$  are the cosh and sinh functions, respectively, and
(\ref{fplus}), (\ref{fminus}) are the familiar sum and difference
relations for these functions. But, of
course, (\ref{gef}) also has many discontinuous solutions
\cite[Chapters 2 and 3]{aczel}.

Conversely, the general solution of (\ref{fminus}) is expressible in
terms of a {\em single} arbitrary exponential function, i.e., it is known
\cite{wilson} that if $f_0$ and $f_1$ satisfy (\ref{fminus}), then
$f_0$ and $f_1$ are the even and odd components of a single  exponential function:
$ f_0(x) = [g(x) + g(-x)]/2,\; f_1(x) = [g(x) - g(-x)]/2. $ In fact,
$g(x) = f_0(x) +f_1(x)$. On the other hand \cite{vietoris}, the
general solution of (\ref{fplus}) depends on {\em two} exponential functions: $
f_0(x) = [g_1(x) + g_2(x)]/2,\; f_1(x) = [g_1(x) - g_2(x)]/2. $ We
could think of $g_1(x)$ as  $e^x$, for
example, and $g_2$ as the identically $0$ function.

The reason that (\ref{fminus}) has fewer solutions than
(\ref{fplus}) is that (\ref{fminus}) implies (\ref{fplus}) but not
vice-versa. To see that (\ref{fminus}) implies (\ref{fplus}), we
first interchange $x$ and $y$ in (\ref{fminus}) to see that $f_0$
must be even and $f_1$ odd. Then replacing $y$ by $-y$ in
(\ref{fminus}) we get (\ref{fplus}). On the other hand, $f_0(x) =
f_1(x) = e^x/2$ satisfies (\ref{fplus}) but not (\ref{fminus}).

 We get
``higher order generalized hyperbolic functions" by considering a
function written as a sum of $n$ components, forming a natural
extension of even and odd components in the case $n =2$.  The
present work discusses the appropriate generalizations of
(\ref{fplus}), (\ref{fminus}) and their solutions in this setting.
Theorem 2.8 gives the systems of equations analogous to
(\ref{fplus}), (\ref{fminus}) and Theorem 3.1 gives the general
solutions of these systems.  Theorem 5.2 discusses the stability of
these systems and makes it clear why the system (\ref{fminus}) is
stable but the system (\ref{fplus}) is not.

The name ``higher order generalized hyperbolic functions" is
intended to convey that the functions discussed here generalise
the hyperbolic functions in two ways.  First, there is the general
order $n$ rather than the special order $2$.  Second, we do not
confine attention to continuous solutions of the relevant
functional equations.

The main results reported in \S 2 and \S 3 are proved in \cite{schw},
essentially by the same methods, but with a less systematic use of
matrix methods.  An attempt to understand \cite{schw} and
\cite{forg} is what led to the present paper.  In \cite{schw} and
\cite{forg}, $n$--tuples of functions were used; here we find it
more convenient to base the development more explicitly on {\em
circulant matrices} formed from these $n$--tuples; this is partially
done in \cite{schw} and \cite{forg}. The question of stability is further considered in \cite{forg1}.

We note that (\ref{fplus}) and (\ref{fminus}) can be written
\begin{equation}
 \left[\begin{array}{cc}
f_0(x+y) & f_1(x+y)\\
f_1(x+y) & f_0(x+y)
\end{array} \right]
=\left[\begin{array}{cc}
f_0(y) & f_1(y)\\
f_1(y) & f_0(y)
\end{array} \right]
\cdot \left[\begin{array}{cc}
f_0(x) & f_1(x)\\
f_1(x) & f_0(x)
\end{array} \right] \label{fplusm}
\end{equation}
and
\begin{equation}
 \left[\begin{array}{cc}
f_0(x-y) & f_1(x-y)\\
f_1(x-y) & f_0(x-y)
\end{array} \right]
=\left[\begin{array}{cc}
f_0(y) & -f_1(y)\\
-f_1(y) & f_0(y)
\end{array} \right]
\cdot \left[\begin{array}{cc}
f_0(x) & f_1(x)\\
f_1(x) & f_0(x)
\end{array} \right], \label{fminusm}
\end{equation}
so our task is to find higher--dimensional analogues of
(\ref{fplusm}) and (\ref{fminusm}).

 Our results are stated
for complex-valued functions of a complex variable, but, as in
\cite{forg}, \cite{forg1} and \cite{schw}, they can be extended to complex-valued
functions on an abelian group $G$ with an automorphism $\sigma$
satisfying $\sigma^n = id_G$.

F. Zorzitto \cite{zorzitto} has considered a far-reaching
generalization, based on group theory, of the results of
\cite{forg} and \cite{schw}.

\section{Preliminary results}
\setcounter{equation}{0}
\subsection{Special matrices}

The $n \times n$ Fourier matrix \cite[p.
32]{davis} ${\cal F}_n$ is given by
\begin{equation}
{\cal F} ={\cal F}_n = \frac{1}{\sqrt{n}}
\left[\begin{array}{ccccc} 1
&1&1 &\dots&1\\
1&\omega^{-1}&\omega^{-2}&\dots&\omega^{-(n-1)}\\
1&\omega^{-2}&\omega^{-4}&\dots& \omega^{-2(n-1)}\\
\dots & \dots &\dots &\dots &\dots \\
1&\omega^{-(n-1)}&\omega^{-2(n-1)}&\dots&\omega^{-(n-1)
^2}
\end{array} \right].
\label{fmat}
\end{equation}
${\cal F}_n$ is a symmetric matrix whose
complex-conjugate is given by
$$
{\cal F}^* ={\cal F}_n^* = \frac{1}{\sqrt{n}}
\left[\begin{array}{ccccc} 1
&1&1 &\dots&1\\
1&\omega&\omega^2&\dots&\omega^{n-1}\\
1&\omega^2&\omega^4&\dots& \omega^{2(n-1)}\\
\dots & \dots &\dots &\dots &\dots \\
1&\omega^{n-1}&\omega^{2(n-1)}&\dots&\omega^{(n-1)
^2}
\end{array} \right].
$$
Sometimes the Fourier matrix is defined as in (\ref{fmat})
but without the factor $1/\sqrt{n}$. The present definition has the
advantage that the matrix ${\cal F}_n $ is unitary:
$
{\cal F}_n^*{\cal F}_n ={\cal F}_n{\cal F}_n^* = I$.
We use the notation \cite{davis} ${\rm circ}({\bf a}^T) = {\rm
circ}(a_0,a_1,\dots,a_{n-1}) $ for the {\em
circulant} matrix:
\begin{equation} A = {\rm circ}({\bf a}^T) =
{\rm circ}(a_0,a_1,\dots,a_{n-1}) = \left[\begin{array}{cccc}
a_0 & a_{1} & \dots & a_{n-1}\\
a_{n-1} & a_{0} & \dots & a_{n-2}\\
\dots & \dots &\dots &\dots \\
a_{1} & a_{2} & \dots & a_0
\end{array} \right].
\label{defa}
\end{equation}

We also use the diagonal matrix
$$
\Omega  = \Omega _n = {\rm diag}[1, \omega,\dots,\omega^{n-1}].
$$
We have \cite[Theorem 3.2.1]{davis}
\begin{equation}
{ \cal F}^{*} \Omega {\cal F} =  \pi,
\label{frm}
\end{equation}
or
\begin{equation}
\Omega =
{ \cal F}\pi {\cal F}^{*},
\label{frm1}
\end{equation}
where $\pi$ is the $n \times n$ permutation matrix
$$ \pi =
\pi_n = {\rm circ}(0,1,\dots,0) = \left[\begin{array}{ccccc}
0 &1&0 &\dots&0\\
0 & 0 & 1 & \dots & 0\\
0 & 0 & 0 & \dots & 0\\
\dots & \dots &\dots &\dots &\dots \\
0 & 0 & 0 & \dots & 1\\
1 & 0 & 0 & \dots & 0
\end{array} \right].
$$
In fact all circulants are polynomials in $\pi$; we have \cite[p.
68]{davis}
\begin{equation}
 {\rm
circ}(a_0,a_1,\dots,a_{n-1})= p_\gamma(\pi),
\label{pif}
\end{equation}
where
\begin{equation}
 p_\gamma(z) = a_0 + a_1z +\cdots + a_{n-1}z^{n-1}. \label{pgz}
\end{equation}

Every circulant matrix can be diagonalized:

{\bf Lemma 2.1}. {\em
If $A$ is a circulant matrix given by (\ref{defa}), then
\begin{equation}
{\cal F} A{\cal F}^* = \sqrt{n}\;{\rm diag}
[{\cal F}^* \;{\bf a}].
\label{l21}
\end{equation}}
Here we  use the usual boldface notation for column vectors:
${\bf a} =
[a_0\; a_1\; \dots \;a_{n-1}]^T
$
and diag[{\bf a}] denotes the square matrix whose diagonal elements
are $a_0, a_1, \dots, a_{n-1}$.

{\em Proof.} This is shown in
(\cite[Theorem 3.2.2]{davis}) in the form
\begin{equation}
{\cal F} A{\cal F}^* = {\rm diag}
(p_\gamma(1), p_\gamma(\omega), \dots,p_\gamma(\omega^{n-1})),
\end{equation}
which is equivalent to (\ref{l21}) on using (\ref{pgz}).

By conjugation and iteration of (\ref{frm}) we get, since
$\pi^T = \pi^{-1}$ and ${\cal F}^T = {\cal F}$,
\begin{equation} {\cal F} \Omega^{-m}{\cal F}^* =
\pi^m, \label{l210}
\end{equation}
for $m= 0,1,\dots$ and,
taking inverses, this also holds for $m = -1,-2,\dots$.
This gives
\begin{equation} {\cal F} [\Omega^{-m}\pi^k \Omega^m]{\cal F}^* =
{\cal F} \Omega^{-m}{\cal F}^*[{\cal F}\pi^k{\cal F}^*]{\cal F} \Omega^m{\cal F}^* =
\pi^m{\cal F} \pi^k {\cal F}^* \pi^{-m}, \label{12110}
\end{equation}
for integers $m$ and nonnegative integers $k$.
Hence, in view of the representation
(\ref{pif}), we get:

{\bf Lemma 2.2.} {\em
For circulant $A$, and $m = 0,1,\dots$,
\begin{equation} {\cal F} [\Omega^{-m}A \Omega^m]{\cal F}^* =
\pi^m{\cal F} A {\cal F}^* \pi^{-m}. \label{1211}
\end{equation}}

\subsection{Components of a function}
Throughout the rest of this paper, we let $n$ be a fixed integer
$\ge 2$.
Following J.~Schwaiger \cite{schw}, we make the following
definition:

{\bf Definition 2.3.} {\em  A function $h: C
\to C$ is
{\em of type}
$j\;(j=0,\dots,n-1$), if $h(\omega x) = \omega^jh(x)$ where $\omega = \omega_n =
e^{2\pi i/n}$.}

  In case $n = 2$, the type $0$ and
type
$1$ functions are the even and odd functions, respectively.
In the cases $n=3$ and $n=4$, examples of functions of types $0,1,2$ and
$0,1,2,3$, respectively, are given in (\ref{inelegant}) and
(\ref{elegant}). It is easy to generate
examples in the case of analytic functions: if $n=3$
and
$f(z) = \sum_{k=0}^\infty a_kz^k,$
in a disk centred at $0$, then
$f_j(z) = \sum_{k=0}^\infty a_{3k+j}z^{3k+j}$
is of type $j,\;j=0,1,2$.   But we have no need of analyticity in
the discussion.

{\bf Lemma 2.4}. {\em
 Every
function $f:C \to C$ can be
expressed uniquely as a sum of functions $f_j$ of type
$j,\;\;j=0,\dots,n-1$, called the {\em components} of $f$. In fact,
\begin{equation}
f = \sum_{j=0}^{n-1}f_j  \label{repn}
\end{equation}
 where
\begin{equation}\left[
\begin{array}{c} f_0(x)\\f_1(x)\\ \vdots \\ f_{n-1}(x)
\end{array} \right]
= \frac{1}{\sqrt{n}}  {\cal F} \left[
\begin{array}{c} f(x)\\f(\omega x)\\ \vdots \\ f(\omega^{
n-1}x)
\end{array} \right]  = \frac{1}{\sqrt{n}} {\cal F}^*\left[
\begin{array}{c} f(x)\\f(\omega^{-1}x)\\ \vdots \\ f(\omega^{-
n+1}x)
\end{array} \right].
 \label{44}
\end{equation}}

{\em Proof.}
Equation (\ref{repn})  with $f_j$ of type $j,\;j=0,\dots,n-1$, implies that
\begin{equation}
\left[
\begin{array}{c} f(x)\\f(\omega x)\\ \vdots \\ f(\omega^{
n-1}x)
\end{array} \right]
= {\sqrt{n}}  {\cal F}^*\left[
\begin{array}{c} f_0(x)\\f_1(x)\\ \vdots \\ f_{n-1}(x)
\end{array} \right]
 \label{444}
\end{equation}
which is obtained from (\ref{repn}) by
replacing $x$ by $x,\omega x, \dots,
\omega^{n-1}x$ in order. Also, of course, (\ref{444}) implies (\ref{repn}).
But (\ref{444})
 obviously has a unique solution for the $f_i$ given
by the first equality in
(\ref{44}), and each of the $f_i$ is of type $i$.

The second equation in (\ref{44}) is obtained similarly.

{\bf Definition 2.5.} {\em The
circulant matrix
\begin{equation}
F(x)= {\rm circ}(f_0,f_1,\dots,f_{n-1}) = \left[\begin{array}{cccc}
f_0 & f_{1} & \dots & f_{n-1}\\
f_{n-1} & f_{0} & \dots & f_{n-2}\\
\dots & \dots &\dots &\dots \\
f_{1} & f_{2} & \dots & f_0
\end{array} \right],
\label{deff}
\end{equation}
whose first
row is
formed from the components of $f$ in order will be called
{\em the circulant matrix corresponding to $f$}.}

{\bf Lemma 2.6.}
{\em If $f:C \to C$ is any function, the   corresponding
circulant matrix function is given by
\begin{equation}
F(x)  = {\cal F}^*{\rm diag}\;[f(x),f(\omega
x),\dots, f(\omega^{n-1}x)]{\cal F}.
\label{l22}
\end{equation}}
{\em Proof.} This comes from Lemma 2.1.

\subsection{Exponential functions}

{\bf Lemma 2.7.}
{\em If
$f$ is an exponential function, the corresponding circulant matrix function $F$ satisfies
$F(x + y) = F(y) F(x).$}

{\em Proof}.
This is a fairly immediate consequence of Lemma 2.6.
${\cal F} F(x+y) {\cal F}^*$ is a diagonal matrix
with entries $f(\omega^{i}(x+y))=f(\omega^{i}x) f(\omega^{
i}y)
,\;i=0,\dots,n-1.$  But $f(\omega^{i}x), f(\omega^{i}y), i
=0,\dots,n-1$, are the respective entries of the diagonal matrices
${\cal F} F(x) {\cal F}^*$
and ${\cal F} F(y) {\cal F}^*$.

{\bf Lemma 2.8.} {\em
For any function
$f$, the corresponding circulant matrix function $F$ satisfies the
equation
$
F(\omega y) = \Omega^{-1} F(y) \Omega.
$
}

{\em Proof.}
This is a fairly immediate consequence of Definition 2.3.

Successive application of Lemmas 2.7 and 2.8 give

{\bf Theorem 2.9.} {\em If
$f$ is an exponential function, the corresponding circulant matrix function $F$ satisfies
\begin{equation}
F(x + \omega^m y) = \Omega^{-m}F(y)\Omega^{m} F(x),
\label{444b}
\end{equation}}
for $m=0,1,\dots,n-1$.

Theorem 2.9 is a matrix form of \cite[equation($4_t$)]{schw}, which
is written in the notation
$$ f_j(x+ \omega^ty) =
\sum_{l=0}^j \omega^{(j-l)t} f_l(x) f_{j-l}(y)
+ \sum_{l=j+1}^{n-l} \omega^{(n+j-l)t} f_l(x) f_{n+j-l}(y), \;t=0,1,\dots,n-1.
$$
In case $n=2$, (\ref{444b}) reduces to (\ref{fplusm}) for $m=0$ and
to (\ref{fminusm}) for $m=1$.

\section{A converse result}
\setcounter{equation}{0}
We now consider equation (\ref{444b}) where
$F(x)$ is circulant, but no assumption is made about its entries $f_i$.
It will come as no surprise that if (\ref{444b}) holds for all
$m=0,1,\dots,n-1$, then the $f_i$ are the components of an exponential function and
$F(x) = {\rm circ}(f_0,f_1,\dots,f_{n-1})$.   But what if we
have (\ref{444b}) for only a single value of $m$?
 Then it is clear that (\ref{444b})
holds also for integer multiples of this value.  Thus if g.c.d.$(n,m) = 1$, and (\ref{444b}) holds for a fixed $m \in \{0,1,\dots,n-1\}$, then it holds for all $m \in \{0,1,\dots,n-1\}$.
  On the other hand, if
(\ref{444b}) holds for a fixed
$m \in \{0,1,\dots,n-1\}$  with
 g.c.d. $(n,m) =d > 1$, then we can assert only that
it holds for values of $m$ which are multiples of $d$, modulo $n$, and
 it turns out that the general solution of (\ref{444b}) depends
on $d$ arbitrary exponential functions:

{\bf Theorem 3.1.} {\em Let a circulant matrix function $F(x)$ satisfy
\begin{equation}
F(x + \omega^m y) = \Omega^{-m}F(y)\Omega^{m} F(x),\label{key1}
\end{equation}
for a {\bf fixed} value of $m  \in \{0,1,\dots,n-1\}$, and let $ d =
{\rm g.c.d.}(n,m)$. Then $F(x) = {\rm circ}({\bf f}(x)^T) $ where
\begin{equation}
{\bf f}(x) = \frac{1}{\sqrt{n}}{\cal F} \left[
\begin{array}{c}{\bf h}(x)\\
 {\bf h}(\omega^{m}x)\\
 \dots \\
{\bf h}(\omega^{m(n/d-1)}x)
\end{array}
\right], \label{3.2}
\end{equation}
${\bf h}(x)$ being an arbitrary $d$--tuple of exponential functions.
In particular, when $d = 1$,
\begin{equation}
{\bf f}(x) = \frac{1}{\sqrt{n}}{\cal F}
\left[
\begin{array}{c}
{g}(x)\\ {g}(\omega x)\\
\dots\\ {g}(\omega^{n-1}x)
\end{array}
\right],  \label{3.3}
\end{equation}
where $g$ is an arbitrary exponential function.}

To prove Theorem 3.1, we use the ideas of \cite{schw}. We
write (\ref{key1}) in the form
\begin{equation}
  {\cal F}F(x + \omega^m y){\cal F}^* = [{\cal F}\Omega^{-m}F(y)\Omega^{m}
  {\cal F}^*][{\cal F} F(x) {\cal F}^*],\label{key11}
\end{equation}
and using Lemmas 2.1 and 2.2, this is
\begin{equation}
G(x+\omega^my) =  \pi^m G(y)\pi^{-m} G(x), \label{key110}
\end{equation}
where $G(x)$ is the diagonal matrix
$$ G(y) = \sqrt{n}\;[{\rm diag}
{\cal F}^*\;{\bf f}(y)].
$$
Using  the notation
\begin{equation}
G_m(x) = \pi^{m}G(x) \pi^{-m}, \label{gandpi}
\end{equation}
we have
\begin{equation}
G_m(x) = {\rm diag}[g_{m}(x), g_{m+1}(x), \dots, g_{m+n-1}(x)],
\label{diag}
\end{equation}
where suffixes are taken modulo $n$, since
$$ \pi
{\rm diag}[\lambda _1,\dots\lambda_n]\pi^{-1} =
{\rm diag}[\lambda _2,\dots\lambda_n,\lambda_1],$$
so the effect of the
transformation $ G(x) \rightarrow \pi^{m}G(x)\pi^{-m}$ is to replace
$g_j$ by $g_{j+m}$.

The easiest case is $m=0$ and hence $d={\rm g.c.d.}(m,n) = n$. In this
case (\ref{key110}) is
 \begin{equation}
G(x+y) =  G(x) G(y), \label{key11n}
\end{equation}
so, since $G$ is diagonal, its entries must be exponential functions. In case
$m=0$ this is all we can say about the $g$'s; they are independent.

The other extreme case is where $d = 1$.
In this case the set of multiples
$\{km,\;k=0,1,\dots,n-1\}$ modulo $n$ is precisely the set $\{1,2,\dots,n\}$ so
(\ref{key110}) can be written as a system of $n$ equations
\begin{equation}
\begin{array}{c}
g_0(x+\omega^m y)  =  g_{m}(y) g_0(x), \\
g_m(x+\omega^m y)  =  g_{2m}(y) g_m(x), \\
\cdots \\
g_{(n-1)m}(x+\omega^m y)  =  g_{0}(y) g_{(n-1)m}(x).
\end{array} \label{3.7}
\end{equation}
There are two possibilities. One is that $g_0(0) = 0$; hence from (\ref{3.7}), all of the $g_i$
are identically $0$ and the assertion of the Theorem
holds trivially. The other possibility is that $g_0(0) \ne 0$.
Putting $x=y=0$ in the first of equations (\ref{3.7}) then gives
$g_{m}(0) = 1$. Then the remaining equations, with $y=0$,  give
$g_{2m}(0) = g_{3m}(0) = \cdots = g_{0}(0) = 1$ or
$g_{2}(0) = g_{3}(0) = \cdots = g_{0}(0) = 1$. Then putting $x = 0$ in (\ref{3.7}) we see that
the $g_k$ are given
by the formula
\begin{equation}
g_{km}(y) = g_0(\omega^{km}y),\; k=0,1,\dots,n-1. \label{310}
\end{equation}
The second of these equations gives
\begin{equation}
g_{m}(y) = g_0(\omega^m y),
\end{equation}
which, together with the first of equations (\ref{3.7}) shows that
$g_0$ is an exponential function.
Now (\ref{310}) is equivalent to
\begin{equation}
g_k(y) = g_0(\omega^{k}y),\; k=0,1,\dots,n-1.
\end{equation}
 Thus $G$ is determined by the single exponential function $g_0$ and ${\bf f}$ is determined by (\ref{3.3}). This completes the proof in the case $m=1$.

In case $d = {\rm g.c.d.}(n,m) > 1$,
the equations (\ref{3.7}) are replaced by $d$ sets of $n/d$ equations corresponding
to the cosets of the additive
 subgroup of $\{0,1,\dots,n-1\}$ generated by $\{0,1,\dots,d-1\}$.
The $r$th such set of equations (where $r$ is one of the numbers $0,1,\dots,d-1$),
reads
\begin{equation}
\begin{array}{c}
g_r(x+\omega^m y)  = g_{r+m}(y) g_r(x),\\
g_{r+m}(x+\omega^m y) =  g_{r+ 2m}(y) g_{r+m}(x),\\
\cdots \\
g_{r+(n/d-1)m}(x+\omega^m y)  =  g_{r+(n/d)m}(y) g_{r+(n/d-1)m}(x).
\end{array}
\end{equation}
Each such set of equations can be dealt with as in the case $m=1$.
That is either all of the $g_{r+ km}\,(k=0,1,\dots,n/d-1) $
are identically $0$ or they
are given by the formula
\begin{equation}
g_{r+ km}(x) = g_r(\omega^{km}x),\; k=0,1,\dots,n/d-1,
\end{equation}
where $g_r$ is an exponential function.
Thus $G$ is determined by the $d$-tuple ${\bf h}(x)^T =
(g_0(x),g_1(x),\dots,g_{d-1}(x))$ of exponential functions and
 we get the main assertion of Theorem 3.1.

\section{ Continuous solutions}
\setcounter{equation}{0}
In the case $n =1$, a continuous solution of $f(x+y) = f(x)f(y)$,
is given by $f(x) = e^{ax}$.
For $n >1$, $e^x$ has components
$$
F^{(1)}_{n,k}(x) = \frac{1}{n}  \sum_{j=0}^{n} \omega^{jk}
e^{\omega^j x},\;k=0,\dots,n-1,
$$
in a minor variation of the
notation of \cite{ung82},
\cite{ung84}. For further references and history of these
functions, see \cite{muldoon}.  It should be noted that the term
``generalized hyperbolic functions" as used in \cite{muldoon},
\cite{ung82} and \cite{ung84}, corresponds to some of the
continuous
solutions of the system of functional equations considered in the
present paper.
The three continuous  generalized hyperbolic functions
of order $3$ are
\begin{equation}
\begin{array}{c}
F_{3,0}^1(x) =
\frac13\left[e^x+2e^{-x/2}\cos\frac{\sqrt{3}x}{2}\right];\\
F_{3,1}^1(x) = \frac13\left[e^x- 2e^{-x/2}\cos
\left(\frac{\sqrt{3}x}{2}+
\frac{\pi}{3}\right)\right];\\
F_{3,2}^1(x) = \frac13\left[e^x- 2e^{-x/2}\cos
\left(\frac{\sqrt{3}x}{2}-\frac{\pi}{3}\right)\right].
\end{array}
\label{inelegant}
\end{equation}

In the case $n=4$, we get the compact formulas
\begin{equation}
\begin{array}{c}
F_{4,0}^{1}(x) = (1/2)(\cosh x + \cos x),\\
F_{4,1}^{1}(x) = (1/2)(\sinh x + \sin x),\\
F_{4,2}^{1}(x) = (1/2)(\cosh x - \cos x),\\
F_{4,3}^{1}(x) = (1/2)(\sinh x - \sin x),
\end{array}
\label{elegant}
\end{equation}
given by Battioni \cite{batt}, a special case of a more
general result
\begin{eqnarray*}  F_{2m,r}^{1}(x) & = &
[F_{m,r}^{1}(x)+F_{m,r}^{-1}(x)]/2, \;\; r=0,1,\dots,m-
1,\\  F_{2m,r+m}^{1}(x) & = & [F_{m,r}^{1}(x)-
F_{m,r}^{-1}(x)]/2, \;\; r=0,1,\dots,m-1,
\end{eqnarray*}
given in a different notation in \cite[(33), p. 216]{bmp}.

The decomposition of various special functions into sums of
functions of order $j,\;j=0,\dots, n-1$ has been considered by Y.
Ben Cheikh; see \cite{bencheikh1} and \cite{bencheikh2}.

\section{Stability}\setcounter{equation}{0}
We use the usual 1-norm for square matrices:
$$ \|A\| =
\|A\|_1 = \max_ {1 \le i \le n}\sum_{j=1}^n |a_{ij}|, $$
If $B$ is a diagonal matrix, we
have
\begin{equation}
\|
AB\| \ge \|A\|
{\rm
min}_{1
\le i \le n}|b_{i,i}|.
\label{ii}
\end{equation}

{\bf Definition 5.1.} (See, e.g., \cite{blz}, \cite{hyers}.) {\em
For fixed $m$, we say that the equation
\begin{equation}
F(x + \omega^m y) = \Omega^{-m}F(y)\Omega^{m} F(x),  \label{444bb}
\end{equation}
for circulant
$F$,
is {\em stable} if the condition
\begin{equation}
\| F(x + \omega^m y) - \Omega^{-m}F(y) \Omega^{m} F(x) \| <
\varepsilon,\;\; {\rm for}\;{\rm some}\; \varepsilon >0 \;{\rm and}\;{\rm all}\;x,y,  \label{4444}
\end{equation}
implies that either $F$ is bounded or $F$ satisfies (\ref{444bb}).}

 It has been noted that in the case $n =2$, the system
(\ref{fminus}) is
stable
\cite{schw91}
but (\ref{fplus}) is not \cite{ger}.  An appropriate
generalization is

{\bf Theorem 5.2.} {\em Let $F$ be a circulant $n \times n$ matrix
function and let $m$ be a fixed integer in the set $\{0,1,\dots, n-1\}$. The equation
\begin{equation}
F(x + \omega^m y) = \Omega^{-m}F(y)\Omega^{m} F(x),  \label{444bbb}
\end{equation}
is stable if and only if g.c.d.$(n,m)=1$}.

{\em Proof.} The idea of the proof is based on the proof of
\cite[Proposition 2]{forg}. Suppose that g.c.d.$(n,m) =1$ and let $F$ be
unbounded and satisfy (\ref{4444}). We need to show that $F$
satisfies (\ref{444bbb}). We
introduce the diagonal matrix
$$
G(x) = {\cal F} F(x) {\cal F}^*
= \sqrt{n}\;{\rm diag}\;[
{\cal F}^*{\bf f}(x)]
$$
as in Lemma 2.1.  Clearly $G$ is also unbounded, since
$ \| {\cal F}\|  = \| {\cal F}^*
\|  = \sqrt{n}. $
These relations, together with Lemma 2.2,
show that (\ref{4444}) is equivalent to
\begin{equation}
\| G(x + \omega^m y) - G_m(y) G(x)\|  \le
 n \varepsilon =\varepsilon_1,
\label{4epsg}
\end{equation}
with the notation of (\ref{gandpi}), (\ref{diag}).
This leads to
$$
\| G(x+ (\omega^my + \omega^mz)) -  G_m(y +z) G(x) \|
\le
\varepsilon_1
$$
and
$$
\| G((x+\omega^mz) + \omega^my) -  G_m(y) G(x+ \omega^m z)
\|  \le \varepsilon_1
$$
so
$$
\|  G_m(y +z) G(x) -  G_m(y) G(x+ \omega^m z)
\|  \le 2\varepsilon_1,
$$
for all $x,y,z$.
Also, from (\ref{4epsg}) and the norm property $\|AB\| \le
\|A\|\|B\|$, we get
$$
\|  G_m(y) G(x+ \omega^m z) -
G_m(y) G_m(z) G(x) \|  \le
\varepsilon_1 \|  G_m(y) \| .
$$
Combining the last two inequalities, we get
$$
\|  G_m(y+z) G(x) -
G_m(y) G_m(z) G(x) \|  \le  \varepsilon_1 (2+
\|  G_m(y) \|).
$$
Since $G$ is diagonal, we may use (\ref{ii}) to get
\begin{eqnarray}
\|  G_m(y +z ) -
G_m(y) G_m(z)\|  {\rm min}_{1 \le i \le n}|g_i(x)| & \le &
\|  [G_m(y +z ) - G_m(y) G_m(z)] G(x)\| \nonumber \\
&\le & \varepsilon_1 (2+
\|  G_m(y) \|). \label{521}
\end{eqnarray}
We see from (\ref{4epsg}) that
$
\| G( \omega^m y) - G_m(y) G(0)\|  \le
\varepsilon_1$;
hence the unboundedness of one diagonal entry in $G$ implies the
unboundedness of {\em all}  diagonal entries in $G$. (Here it is
crucial that g.c.d.$(n,m)=1$; otherwise integers multiples of $m$
will
not cover all $j=0,\dots,n-1$.) This means
that,
given $y$ and $z$, we can choose $x$ such that
$$
{\rm min}_{1 \le i \le n}|g_i(x)| \ge 2+ \|
G_m(y)\|.
$$
Hence, from (\ref{521}),
$
\|  G_m(y +z ) - G_m(y) G_m(z)\|
 \le \varepsilon_1.
$
Since the elements of $G_m$ are the same as those of $G$, in
a different order, we finally get
\begin{equation}
\|  G(y +z ) - G(y) G(z)\|   \le \varepsilon_1.
\label{4444g}
\end{equation}
Thus each element $g_i$ in the diagonal matrix $G(y)$ is
unbounded and satisfies the inequality
$
|g_i(y +z ) - g_i(y) g_i(z)|  \le \varepsilon_1.
$
By  \cite[Theorem 1]{baker}, we must therefore have
$
 g_i(y +z ) - g_i(y) g_i(z) =0,
$
so
\begin{equation}
 G(y +z ) - G(y) G(z) =0  \label{444g}
\end{equation}
and, as in the proof of Theorem 2.9, $F$ is a solution of
(\ref{444bb}).

To prove the ``only if" assertion of the Theorem, we note that if
g.c.d.$(n,m) > 1$, the unboundedness of one
diagonal entry in $G$, does not imply the unboundedness of all
entries there.   So, for example, suppose that $g_i(x) = e^x$, for
$i$ a multiple of $d = {\rm g.c.d.}(n,m) > 1$ and $g_i(x) =
2$,
when $i$ is not
a multiple of $d$.  Then $G$ is an unbounded solution of
(\ref{4444g}), with $\varepsilon_1 = 2$,  which does not satisfy
(\ref{444g}). Thus we can have
an unbounded solution of (\ref{4444}) which does not satisfy
(\ref{444bbb}).
This, of course, is exemplified by the non-stability of the
equations (\ref{fplus}).

{\bf Acknowledgments.} I thank the referees for very
helpful suggestions and for identifying some errors in  earlier
versions of this article.

\end{document}